\documentclass[11pt,letterpaper]{amsart}

\usepackage[T1]{fontenc}
\usepackage{lmodern}
\usepackage{amsmath,amssymb,amsthm,mathtools}
\usepackage[margin=1.12in]{geometry}
\usepackage{microtype}
\usepackage{enumitem}
\usepackage{longtable}
\usepackage[colorlinks=true,linkcolor=blue!45!black,citecolor=blue!45!black,urlcolor=blue!45!black]{hyperref}
\usepackage{xcolor}

\setlist[enumerate]{leftmargin=*,itemsep=0.5em,topsep=0.5em}
\newtheorem{theorem}{Theorem}
\newtheorem{corollary}{Corollary}
\newtheorem*{nashconjecture}{Nash blowup conjecture}
\newtheorem*{openproblem}{Open problem}
\newtheorem{lemma}{Lemma}
\newtheorem{claim}{Claim}
\theoremstyle{remark}
\newtheorem{remark}{Remark}

\newcommand{\NN}{\mathbb Z_{\ge 0}}
\newcommand{\RR}{\mathbb R}
\newcommand{\Aaff}{\mathbb A}
\DeclareMathOperator{\Cone}{Cone}
\DeclareMathOperator{\Spec}{Spec}
\DeclareMathOperator{\Nash}{Nash}
\DeclareMathOperator{\intt}{int}

\title[Isolated Singularities and Iterated Nash Blowups]{Isolated Singularities Not Resolved\\by Iterated Nash Blowups}

\author{Tongyan Chen}
\address{Zhili College, Tsinghua University, Beijing, 100084, P. R. China}
\email{chen-ty23@mails.tsinghua.edu.cn}

\author{Huaiqing Zuo}
\address{Department of Mathematical Sciences, Tsinghua University, Beijing, 100084, P. R. China}
\email{hqzuo@mail.tsinghua.edu.cn}
\date{}
\hypersetup{pdftitle={Isolated Singularities Not Resolved by Iterated Nash Blowups},pdfauthor={Tongyan Chen and Huaiqing Zuo}}

\begin{document}
\begin{abstract}
Castillo, Duarte, Leyton-\'Alvarez, and Liendo asked whether iterating ordinary Nash blowups in characteristic zero or normalized Nash blowups in arbitrary characteristic can fail to resolve isolated singularities~\cite[Section~4.3]{computational}. We answer this question affirmatively with two constructions of isolated singularities, one nonnormal and one normal, both giving counterexamples to normalized Nash iteration in every characteristic.


\end{abstract}
\maketitle
\begingroup
\renewcommand{\thefootnote}{}
\footnotetext{H. Zuo was supported by NSFC Grant (No. 12671056) and BJNSF Grant
	(No. 1252009).}
\endgroup

\section{Introduction}\label{sec:introduction}

Throughout, $k$ is an algebraically closed field. For a $d$-dimensional affine variety $X\subset\Aaff_k^n$, the \emph{Nash blowup} $\Nash(X)$ is the closure in $X\times\operatorname{Gr}(d,n)$ of the graph of the Gauss map $x\mapsto T_xX$ from the smooth locus of $X$ to the Grassmannian of $d$-planes in $k^n$. The natural projection $\nu_X:\Nash(X)\to X$ is a proper birational morphism that is an isomorphism over the smooth locus. This construction is independent of the chosen embedding and glues on affine open subsets, defining the Nash blowup of any variety. The \emph{normalized Nash blowup} is the normalization of $\Nash(X)$, equipped with the induced morphism to $X$.

Nash proposed resolving singularities by iterating the Nash blowup~\cite{nash}, leading to the following conjecture.

\begin{nashconjecture}
Assume $\operatorname{char}(k)=0$. For every variety $X$ over $k$, there exists an integer $m\ge0$ such that the variety obtained from $X$ by $m$ successive ordinary Nash blowups is smooth.
\end{nashconjecture}

Its normalized version made the same assertion in arbitrary characteristic, with normalization after each Nash blowup~\cite{spivakovsky}.

In characteristic zero, Nobile proved that the Nash blowup is an isomorphism precisely when the variety is smooth, and deduced that iterated ordinary Nash blowups resolve curves~\cite{nobile}. Gonz\'alez-Sprinberg introduced the normalized approach~\cite{gonzalez-sprinberg}, and Spivakovsky subsequently proved that iterated normalized Nash blowups resolve complex surfaces~\cite{spivakovsky}. For toric varieties, Gonz\'alez P\'erez and Teissier extended the combinatorial description of the Nash blowup to the nonnormal setting~\cite{gonzalez-teissier}.

In positive characteristic, Nobile exhibited singular curves isomorphic to their Nash blowups, so ordinary iteration can fail~\cite{nobile}. Their normalizations are smooth, so they do not settle the normalized question. Duarte and N\'u\~nez-Betancourt extended Nobile's smoothness criterion to normal varieties~\cite{duarte-nunez}. Positive results include resolution of normal toric surfaces by iterated normalized Nash blowups~\cite{duarte-prime}, smoothness of the ordinary Nash blowup of $2$-generic determinantal varieties~\cite{dalbelo-duarte-ruas}, and further toric families resolved by one normalized Nash blowup in arbitrary characteristic~\cite{characteristic-free}.

Castillo, Duarte, Leyton-\'Alvarez, and Liendo disproved the ordinary Nash blowup conjecture in characteristic zero and its normalized version in every characteristic, in dimensions at least four in their beautiful paper~\cite{counterexamples}. Their four-dimensional examples are the normal affine toric varieties $\Spec k[K_i\cap\mathbb Z^4]$, where $K_i\subset\RR^4$ is generated by the columns of $B_i$, for $i=1,2,3$. The three matrices are those in~\cite[formulas~(4), (6), and~(7)]{counterexamples}:
\begin{equation}\label{eq:recurring-cones}
\begin{gathered}
K_i=\Cone(B_i)\quad(i=1,2,3),\\[4pt]
B_1=\begin{pmatrix}
1&0&0&0&2&1\\
0&1&0&0&3&3\\
0&0&1&0&-2&-1\\
0&0&0&1&-1&-1
\end{pmatrix},\quad
B_2=\begin{pmatrix}
1&0&1&0&0\\
0&1&1&0&2\\
0&0&2&0&2\\
0&0&0&1&1
\end{pmatrix},\quad
B_3=\begin{pmatrix}
1&0&1&0&0\\
0&1&2&0&3\\
0&0&3&0&3\\
0&0&0&1&1
\end{pmatrix}.
\end{gathered}
\end{equation}
A chart means an affine open subset. We write $X\rightsquigarrow Y$ if $Y$ is isomorphic to a chart of $\Nash(X)$, and $X\overset{n}{\rightsquigarrow}Y$ if $Y$ is isomorphic to a chart of the normalization of $\Nash(X)$. A \emph{chart path} is a finite sequence of these choices with a fixed arrow type. A path's \emph{source} is its initial variety. An ordinary arrow with normal target also gives the normalized arrow, since normalization is an isomorphism over that chart.

A singular affine variety $Y$ is called \emph{recurring} under ordinary, respectively normalized, Nash blowups if some positive iterate of the corresponding operation on $Y$ contains a chart isomorphic to $Y$.

For a lattice $M\cong\mathbb Z^d$, write $M_{\RR}=M\otimes_{\mathbb Z}\RR$. A cone $K\subset M_{\RR}$ is \emph{rational polyhedral} if it consists of all nonnegative real linear combinations of finitely many vectors in $M$; it is \emph{full-dimensional} if its real linear span is $M_{\RR}$, and \emph{pointed} if it contains no line through the origin.

For a full-dimensional pointed rational polyhedral cone $K\subset M_{\RR}$, set $S_K=K\cap M$ and $X_K=\Spec k[S_K]$. Here the cone subscript means that we take the semigroup of lattice points in the cone.

\begin{lemma}[{\cite[Section~2]{counterexamples}}]\label{lem:recurring-charts}
The varieties $X_{K_1},X_{K_2},X_{K_3}$ are normal and singular, and have the following recurring chart paths:
\[
\begin{aligned}
X_{K_1}&\rightsquigarrow X_{K_1},
&&\operatorname{char}(k)\notin\{2,3\},\\
X_{K_2}&\overset{n}{\rightsquigarrow}X_{K_2},
&&\operatorname{char}(k)=2,\\
X_{K_3}&\overset{n}{\rightsquigarrow}X_{K'_3}\overset{n}{\rightsquigarrow}X_{K_3},
&&\operatorname{char}(k)=3.
\end{aligned}
\]
The first path also holds with $n$-arrows, since its selected chart is normal. Here $K'_3$ is the cone of the intermediate normal affine chart $X_{K'_3}$ in characteristic $3$, with generators given in~\eqref{eq:connecting-cones-three}. This chart is singular because its normalized Nash blowup contains the singular chart $X_{K_3}$.
\end{lemma}

The ordinary recurrence for $X_{K_1}$ in characteristic $p>3$ follows from the same computation as in characteristic zero: all nonzero replacement determinants in~\cite[formula~(5)]{counterexamples} have absolute value $1$, $2$, or $3$, so the selected chart semigroup is unchanged.

The same authors later found a nonnormal affine variety $Z$ of dimension $3$ with the ordinary chart path~\cite{ordinary3}
\[
Z\rightsquigarrow Z'\rightsquigarrow Z,
\qquad\operatorname{char}(k)=0.
\]
Here $Z'$ denotes the intermediate chart. The varieties $X_{K_1},X_{K_2},X_{K_3}$ and $Z$ have non-isolated singularities.

Thus these examples did not rule out the possibility that the Nash blowup conjectures might hold for varieties with isolated singularities. This led the authors to pose the following problem~\cite[Section~4.3, fourth item]{computational}.
\begin{openproblem}[Castillo, Duarte, Leyton-\'Alvarez, and Liendo]
Can iterating ordinary Nash blowups in characteristic zero, or normalized Nash blowups in arbitrary characteristic, fail to resolve an isolated singularity?
\end{openproblem}
Our results give an affirmative answer in both settings. Our first result realizes a prescribed toric chart while placing all singularities of the source at one point.

\begin{theorem}[Prescribed-chart realization]\label{thm:main}
Let $d\ge1$, let $M\cong\mathbb Z^d$, and let $K\subset M_{\RR}$ be a full-dimensional pointed rational polyhedral cone. There exists an affine variety $X(K)$ of dimension $d$ over $k$, with an isolated singularity and normalization $\Aaff_k^d$, such that:
\begin{enumerate}
\item if $\operatorname{char}(k)=0$, then $X_K$ is a normal chart of the Nash blowup of $X(K)$;
\item if $\operatorname{char}(k)>0$, then $X_K$ is a chart of the normalized Nash blowup of $X(K)$.
\end{enumerate}
\end{theorem}

Applying Theorem~\ref{thm:main} to the cones $K_1,K_2,K_3$ defined in~\eqref{eq:recurring-cones} supplies the first arrow in each of the following paths, attaching a source with an isolated singularity to the corresponding recurrence in Lemma~\ref{lem:recurring-charts}:
\[
\begin{aligned}
X(K_1)&\rightsquigarrow X_{K_1}\rightsquigarrow X_{K_1},
&&\operatorname{char}(k)=0,\\
X(K_1)&\overset{n}{\rightsquigarrow}X_{K_1}\overset{n}{\rightsquigarrow}X_{K_1},
&&\operatorname{char}(k)\notin\{2,3\},\\
X(K_2)&\overset{n}{\rightsquigarrow}X_{K_2}\overset{n}{\rightsquigarrow}X_{K_2},
&&\operatorname{char}(k)=2,\\
X(K_3)&\overset{n}{\rightsquigarrow}X_{K_3}\overset{n}{\rightsquigarrow}X_{K'_3}\overset{n}{\rightsquigarrow}X_{K_3},
&&\operatorname{char}(k)=3.
\end{aligned}
\]
In characteristic zero, both selected ordinary charts are normal, so the first path also holds with $n$-arrows. Since Nash blowup and normalization commute with restriction to open subsets, these recurrences produce singular charts at arbitrarily late stages. Both operations are isomorphisms on smooth varieties, so no iterate can be smooth.

We therefore obtain the following corollary.

\begin{corollary}[Nonnormal isolated counterexamples]\label{cor:nonnormal-isolated}
Let $k$ be an algebraically closed field. There exists a nonnormal affine variety $X$ of dimension $4$ over $k$, with an isolated singularity and normalization $\Aaff_k^4$, such that every iterate of its normalized Nash blowup is singular. If $\operatorname{char}(k)=0$, every iterate of its ordinary Nash blowup is also singular.
\end{corollary}

The normalized Nash blowup is the normalization of $\Nash(X)$; it can be singular even when the normalization of $X$ is smooth.

Corollary~\ref{cor:nonnormal-isolated} answers the open problem of Castillo, Duarte, Leyton-\'Alvarez, and Liendo affirmatively with nonnormal counterexamples. Theorem~\ref{thm:main} gives a uniform construction in arbitrary dimension, with potential applications beyond these examples. Our second result gives normal isolated counterexamples to normalized Nash iteration in every characteristic.

\begin{samepage}
\begin{theorem}[Normal isolated counterexample]\label{thm:normal-isolated-cycle}
Let $k$ be an algebraically closed field. There exists a normal Gorenstein affine toric variety $X$ of dimension $4$ over $k$, with an isolated singularity, such that every iterate of its normalized Nash blowup is singular. If $\operatorname{char}(k)\notin\{2,3\}$, every iterate of its ordinary Nash blowup is also singular.
\end{theorem}
\end{samepage}

The proof uses the following paths. We write $C_i^{(a)}$ for the source and intermediate cones specified in~\eqref{eq:source-cone-one}--\eqref{eq:intermediate-cone-two}, and set $X_i^{(a)}=X_{C_i^{(a)}}$, where $a\in\{1,2,3\}$ labels the recurrence containing $X_{K_a}$ and $i$ counts steps from the source, starting at $0$. The recurring cones retain their $K$-notation. The first and third paths share their initial cone models: $C_i^{(3)}=C_i^{(1)}$ for $i=0,1$. The sources and connecting charts are described in Section~\ref{sec:normal-cycle}:
\begin{equation}\label{eq:normal-chart-paths}
\begin{aligned}
X_0^{(1)}&\rightsquigarrow X_1^{(1)}\rightsquigarrow X_2^{(1)}\rightsquigarrow X_3^{(1)}
\rightsquigarrow X_{K_1}\rightsquigarrow X_{K_1},
&&\operatorname{char}(k)\notin\{2,3\},\\
X_0^{(2)}&\overset{n}{\rightsquigarrow}X_1^{(2)}\overset{n}{\rightsquigarrow}X_{K_2}\overset{n}{\rightsquigarrow}X_{K_2},
&&\operatorname{char}(k)=2,\\
X_0^{(3)}&\overset{n}{\rightsquigarrow}X_1^{(3)}\overset{n}{\rightsquigarrow}X_2^{(3)}
\overset{n}{\rightsquigarrow}X_{K'_3}\overset{n}{\rightsquigarrow}X_{K_3}
\overset{n}{\rightsquigarrow}X_{K'_3},
&&\operatorname{char}(k)=3.
\end{aligned}
\end{equation}
The selected ordinary charts in the first path are normal, so this path also holds with $n$-arrows. Normalization is required for some charts in the other two paths. The recurring arrows in all three paths are supplied by Lemma~\ref{lem:recurring-charts}.

We begin with the toric chart formula in Section~\ref{sec:toric}. Sections~\ref{sec:construction} and~\ref{sec:verification} construct and verify the prescribed chart and apply the construction to $K_1,K_2,K_3$. Section~\ref{sec:normal-cycle} constructs normal isolated examples and verifies their chart paths in every characteristic. Appendix~\ref{app:ordinary-refinement} records a positive-characteristic construction refinement.

\begin{samepage}
\begin{center}
\bf{Acknowledgements}
\end{center}

 AI tools (ChatGPT and DeepSeek) were used during the preparation of the paper as technical assistants to help with checking calculations and arguments, and improving the exposition. All mathematical statements, proofs, and conclusions were verified by the authors.
\end{samepage}

\section{Toric Nash charts}\label{sec:toric}

Let $M\cong\mathbb Z^d$ be a lattice. An affine semigroup $P\subset M$ is a finitely generated additive submonoid containing $0$; we assume that $P$ generates $M$ as a group. Its cone is
\[
\Cone(P)=\left\{\sum_{u\in P}\lambda_u u\ \middle|\ \lambda_u\in\RR_{\ge0}\text{ and only finitely many }\lambda_u\ne0\right\}\subset M_{\RR}.
\]
We call $P$ pointed if $\Cone(P)$ is pointed. Its saturation is
\[
\overline P=\Cone(P)\cap M;
\]
we call $P$ saturated if $P=\overline P$.

The semigroup algebra $k[P]$ has a $k$-basis $\{\chi^u\mid u\in P\}$, with multiplication $\chi^u\chi^v=\chi^{u+v}$. The affine toric variety defined by $P$ is
\[
X_P=\Spec k[P].
\]
It has dimension $d$, and it is normal precisely when $P$ is saturated. For a full-dimensional pointed rational polyhedral cone $K\subset M_{\RR}$, the notation $X_K$ is shorthand for $X_{S_K}$, where $S_K=K\cap M$. Since $S_K$ is saturated, $X_K$ is normal.

Assume henceforth that $P$ is pointed, let $H$ be its Hilbert basis, that is, the set of its nonzero indecomposable elements, and write $M_k=M\otimes_{\mathbb Z}k$. Call a subset $A=\{a_1,\ldots,a_d\}\subset H$ a \emph{$k$-basis} if its image is a basis of $M_k$. Equivalently, $\det(a_1,\ldots,a_d)\ne0$ in $k$; this condition may depend on the characteristic of $k$.

For each $k$-basis $A\subset H$, define
\[
\Gamma_A=\left\langle H\cup\left\{h-a\ \middle|\ a\in A,\ h\in H\setminus A,\ (A\setminus\{a\})\cup\{h\}\text{ is a $k$-basis}\right\}\right\rangle_{\NN}.
\]

\begin{lemma}\label{lem:chart}
$\Spec k[\Gamma_A]$ is an affine chart of the Nash blowup of $X_P$ whenever $\Gamma_A$ is pointed. These charts cover the Nash blowup as $A$ ranges over the $k$-bases for which $\Gamma_A$ is pointed.
\end{lemma}
\begin{proof}
See~\cite[Propositions~32 and~60]{gonzalez-teissier} in characteristic zero and~\cite[Theorem~1.9]{duarte-prime} in positive characteristic. A unified combinatorial description is given in~\cite[Section~1]{counterexamples}.
\end{proof}

\begin{remark}\label{rem:generators}
The Hilbert basis in the definition of $\Gamma_A$ may be replaced by any finite generating set $G$ of $P$ containing $A$: the resulting semigroup is unchanged. Indeed, if a generator $g\in G\setminus A$ gives a valid replacement of $a\in A$, express $g$ as a sum of Hilbert-basis elements. Multilinearity of the determinant implies that at least one summand $h$ gives a nonzero replacement determinant. If $h=a$, then $g-a\in P$. Otherwise $h\in H\setminus A$ and
\[
g-a=(h-a)+(g-h)\in\Gamma_A,
\]
since $g-h\in P$. The reverse inclusion follows from $H\subset G$, since every generating set of a pointed semigroup contains its nonzero indecomposable elements.
\end{remark}

\section{A source for a prescribed toric chart}\label{sec:construction}

We first build a semigroup whose associated affine variety has a unique singular point at the closed orbit. The distinguished generators produce the desired chart; adding lattice points of sufficiently large degree makes the normalization an isomorphism away from that orbit.

Fix a full-dimensional pointed rational polyhedral cone $K\subset M_{\RR}$, and write $S_K=K\cap M$. Let $H=\{h_1,\ldots,h_s\}$ be the Hilbert basis of $S_K$.

Choose a primitive lattice point $v_1\in\intt(K)$, and extend it to a lattice basis $v_1,v_2',\ldots,v_d'$ of $M$. For each $i=2,\ldots,d$, choose a sufficiently large $c_i\in\NN$ and set $v_i=v_i'+c_i v_1$. Then $v_i\in\intt(K)$, and $v_1,\ldots,v_d$ remains a lattice basis of $M$. Write each $m\in M$ uniquely as $m=\sum_{i=1}^d m^i v_i$, with $m^i\in\mathbb Z$, and define $l(m)=\sum_{i=1}^d m^i$. Then $l$ is $\mathbb Z$-linear.

Let $K_0=\Cone(v_1,\ldots,v_d)\subset K$. Given $u_0\in K_0\cap M$, define $u_j$, for $j=1,\ldots,s$, by
\begin{equation}\label{eq:shift}
u_j=\begin{cases}
u_0+h_j,&\operatorname{char}(k)=0,\\
u_0+p h_j,&\operatorname{char}(k)=p>0.
\end{cases}
\end{equation}
Since $u_0,h_j\in K\cap M$, each $u_j$ is also a lattice point of $K$.

\begin{claim}
We can choose $u_0\in K_0\cap M$ so that
\[
u_j\in K_0\quad(j=0,\ldots,s),\qquad
u_0-v_i\in K_0\quad(i=1,\ldots,d-1),\qquad
u_j^d\ge2\quad(j=0,\ldots,s).
\]
If $\operatorname{char}(k)=p>0$, we can also require
\begin{equation}\label{eq:congruence}
u_0\equiv v_d\pmod{pM}.
\end{equation}
\end{claim}
\begin{proof}
In the basis $v_1,\ldots,v_d$, the cone $K_0$ corresponds to $\RR_{\ge0}^d$, and $K_0\cap M$ corresponds to $\mathbb Z_{\ge0}^d$. Thus $u_0-v_i\in K_0$ for $i<d$ whenever $u_0^i\ge1$. The conditions $u_j\in K_0$ and $u_j^d\ge2$ amount to finitely many lower bounds on the coordinates of $u_0$: for $j\ge1$, $u_j^i=u_0^i+h_j^i$ in characteristic zero and $u_j^i=u_0^i+p h_j^i$ in characteristic $p>0$. In positive characteristic, choose each $u_0^i$, $i<d$, to be a sufficiently large multiple of $p$, and choose $u_0^d$ sufficiently large with $u_0^d\equiv1\pmod p$. This also ensures $u_0\equiv v_d\pmod{pM}$. In characteristic zero, choose all coordinates of $u_0$ sufficiently large. In either case, these choices meet all the finitely many coordinate lower bounds, so the claim follows.
\end{proof}

\begin{remark}
In characteristic zero, Theorem~\ref{thm:main} realizes $X_K$ as a normal chart of the ordinary Nash blowup of $X(K)$. In positive characteristic, the selected exchanges produce $ph_j$ rather than $h_j$ in~\eqref{eq:shift}, so the construction guarantees $X_K$ as a chart only after normalization. This conclusion is sufficient to answer the open problem of Castillo, Duarte, Leyton-\'Alvarez, and Liendo in positive characteristic, since that part concerns normalized Nash blowups; see Corollary~\ref{cor:nonnormal-isolated}. A refinement for $K_1$ and $K_3$ in suitable positive characteristics is discussed in Appendix~\ref{app:ordinary-refinement}.
\end{remark}

Set
\[
A=\{v_1,\ldots,v_{d-1},u_0\}.
\]
Choose a finite family of nonzero forms $f_1,\ldots,f_m\in\operatorname{Hom}(M,\mathbb Z)$ defining $K$:
\[
K=\bigcap_{\alpha=1}^m\{x\in M_{\RR}:f_\alpha(x)\ge0\}.
\]
For each $\alpha=1,\ldots,m$, set
\[
r_\alpha=\frac{\max_{a\in A}f_\alpha(a)}{\min_{1\le i\le d}f_\alpha(v_i)}.
\]
The denominator is positive because each $v_i\in\intt(K)$ and each $f_\alpha$ is a nonzero form defining $K$. Choose an integer $N_0$ sufficiently large that
\[
N_0>\max_{1\le\alpha\le m}r_\alpha.
\]
We also require $N_0>1$ and $N_0>l(u_0)$.

Finally, let
\[
Q=\NN\{v_1,\ldots,v_{d-1},u_0,\ldots,u_s\},\qquad
T_{N_0}=\{u\in K_0\cap M:l(u)\ge N_0\},\qquad P=Q\cup T_{N_0}.
\]
The claim implies $Q\subset K_0$ and $l(q)\ge0$ for every $q\in Q$. Also $T_{N_0}+T_{N_0}\subset T_{N_0}$, and $Q+T_{N_0}\subset T_{N_0}$. Hence $P$ is closed under addition and contains $0$, so it is an additive semigroup.

\section{Proof of the prescribed-chart theorem}\label{sec:verification}

We first check the source variety, then identify the selected Nash chart.

\subsection{An isolated singularity with smooth normalization}

We first verify the following properties of $P$.
\begin{enumerate}
\item \textbf{Finite generation.} In the basis $v_1,\ldots,v_d$, $T_{N_0}$ consists of the vectors in $\mathbb Z_{\ge0}^d$ with coordinate sum at least $N_0$. Let
\[
F_{N_0}=\{t\in K_0\cap M:N_0\le l(t)<2N_0\}
\]
and set
\[
G=\{v_1,\ldots,v_{d-1},u_0,\ldots,u_s\}\cup F_{N_0}.
\]
Every element of $T_{N_0}$ is a sum of elements of the finite set $F_{N_0}$: repeatedly split off a nonnegative subvector of coordinate sum $N_0$ until the remainder has degree between $N_0$ and $2N_0-1$. Thus $G$ generates $P$.

\item \textbf{Pointedness, dimension, and normalization.} Since $P\subseteq K_0$, it is pointed. For each $i$, both $N_0v_i$ and $(N_0+1)v_i$ lie in $P$, so the group generated by $P$ contains their difference $v_i$ and is all of $M$. Also $N_0v_i\in P$ for every $i$, hence $\Cone(P)=K_0$. Therefore $P$ has rank $d$ and saturation
\[
\overline P=\Cone(P)\cap M=K_0\cap M=\NN\{v_1,\ldots,v_d\}.
\]
Because $v_1,\ldots,v_d$ is a lattice basis, $k[\overline P]\cong k[x_1,\ldots,x_d]$. Set $X(K)=X_P$. This variety has dimension $d$ and normalization $\Aaff_k^d$.

\item \textbf{The singular locus.} Let $\overline{\mathfrak m}=(\chi^{v_1},\ldots,\chi^{v_d})$ in $k[\overline P]$. Since every lattice point of $K_0$ of $l$-degree at least $N_0$ belongs to $T_{N_0}\subset P$, the ideal $\overline{\mathfrak m}^{\,N_0}$ lies in $k[P]$ and in the conductor of $k[P]\subset k[\overline P]$. Its radical in $k[P]$ is the maximal ideal $\mathfrak m_P$ of the origin: every nonconstant monomial of $k[P]$ has positive $l$-degree, so some power belongs to $\overline{\mathfrak m}^{\,N_0}$. Thus the quotient $k[\overline P]/k[P]$, which is annihilated by the conductor, vanishes after localization at every prime other than $\mathfrak m_P$. The normalization is therefore an isomorphism away from the origin.

On the other hand, $v_d\notin P$: it has $l$-degree $1<N_0$ and $d$-th coordinate $1$, whereas the generators of $Q$ have $d$-th coordinate $0$ or at least $2$. Hence $\chi^{v_d}\in k[\overline P]\setminus k[P]$ is integral over $k[P]$. It remains absent after localization at $\mathfrak m_P$, since multiplying it by an element with nonzero constant term leaves a nonzero $\chi^{v_d}$-term. Thus the origin is nonnormal, hence singular, and it is the unique singular point.
\end{enumerate}

\subsection{Recovering the prescribed cone}

We then verify the properties needed for $A$ to index the desired chart.
\begin{enumerate}
\item \textbf{Hilbert-basis membership.} Each $v_i$, $i<d$, is indecomposable: $l(v_i)=1$, while every nonzero element of $P$ has positive $l$-degree. Hence these $v_i$ belong to the Hilbert basis of $P$.

We also show that $u_0$ is indecomposable. Suppose $u_0=x+y$ with nonzero $x,y\in P$. Since $l(u_0)<N_0$, neither summand lies in $T_{N_0}$, so $x,y\in Q$ and can be written using $v_1,\ldots,v_{d-1},u_0,\ldots,u_s$. If this sum contains some $u_j$ with $j\ge1$, then $u_j=u_0+h_j$ in characteristic zero and $u_j=u_0+p h_j$ in characteristic $p>0$. Since $h_j\ne0$ and $K$ is pointed, some defining form $f_\alpha$ satisfies $f_\alpha(h_j)>0$. Thus $f_\alpha(u_j)>f_\alpha(u_0)$; all other summands lie in $K$, so their $f_\alpha$-values are nonnegative. The sum therefore cannot equal $u_0$. If the sum contains $u_0$ and another nonzero generator, that generator lies in $\intt(K)$, so every defining form is positive on it; again the sum cannot equal $u_0$. Finally, a sum using only $v_1,\ldots,v_{d-1}$ has $d$-th coordinate zero, unlike $u_0$, for which $u_0^d\ge2$. Therefore $u_0$ is indecomposable and belongs to the Hilbert basis of $P$.

\item \textbf{Basis condition.} With the fixed ordering $v_1,\ldots,v_{d-1},u_0$, the determinant of the vectors in $A$ in the basis $v_1,\ldots,v_d$ is $u_0^d$. It is nonzero in characteristic zero. In characteristic $p>0$, the earlier choice $u_0\equiv v_d\pmod{pM}$ gives $u_0^d\equiv1\pmod p$, so the determinant is nonzero modulo $p$. Therefore $A$ is a basis of $M\otimes_{\mathbb Z}k$.

\item \textbf{The chart semigroup lies in $K$.} By the finite-generation argument in Section~\ref{sec:verification}, $G$ is a finite generating set of $P$. By Remark~\ref{rem:generators}, we may use it, rather than the Hilbert basis, to compute the chart semigroup:
\[
\Gamma_A=\left\langle G\cup\left\{g-a\ \middle|\ a\in A,\ g\in G\setminus A,\ (A\setminus\{a\})\cup\{g\}\text{ is a $k$-basis}\right\}\right\rangle_{\NN}.
\]
First, every element of $G$ lies in $K$: this holds for the $v_i$ and $u_j$ by construction, and for $F_{N_0}$ since $K_0\subset K$. For the exchange differences, first consider the distinguished generators $u_j$, $j\ge1$. For $a=u_0$, $u_j-a=h_j\in K$ in characteristic zero and $u_j-a=ph_j\in K$ in characteristic $p>0$. For $a=v_i$, the difference lies in $K$, since
\[
u_j-a=\begin{cases}
(u_0-v_i)+h_j,&\operatorname{char}(k)=0,\\
(u_0-v_i)+p h_j,&\operatorname{char}(k)=p>0.
\end{cases}
\]
Here $u_0-v_i\in K_0\subset K$ by the choice of $u_0$. Finally, let $t\in F_{N_0}$ and $a\in A$. Write $t=\sum_i t^i v_i$. For each $\alpha$, we have
\begin{align*}
f_\alpha(t-a)
&=\sum_i t^i f_\alpha(v_i)-f_\alpha(a)\\
&\ge l(t)\min_i f_\alpha(v_i)-\max_{a'\in A}f_\alpha(a')\\
&\ge N_0\min_i f_\alpha(v_i)-\max_{a'\in A}f_\alpha(a')>0,
\end{align*}
by the choice $N_0>r_\alpha$. Thus every valid exchange difference lies in $K$, and consequently $\Gamma_A\subseteq K$.

\item \textbf{Identification of the chart.} For each $j=1,\ldots,s$, replacing $u_0$ by $u_j$ in $A$ gives determinant $u_j^d$. This is positive in characteristic zero and congruent to $1$ modulo $p$ in characteristic $p>0$, so each replacement is valid. Hence $\Gamma_A$ contains
\[
u_j-u_0=\begin{cases}
h_j,&\operatorname{char}(k)=0,\\
ph_j,&\operatorname{char}(k)=p>0.
\end{cases}
\]
Since $P\subseteq\Gamma_A$ and $P$ generates $M$ as a group, so does $\Gamma_A$. In characteristic zero, these $h_j$ generate $S_K$. Since $\Gamma_A\subseteq K\cap M=S_K$, it follows that $\Gamma_A=S_K$. In characteristic $p>0$, the vectors $ph_j$ generate the cone $K$, so $\Cone(\Gamma_A)=K$ and $\overline{\Gamma_A}=K\cap M=S_K$. The inclusion $\Gamma_A\subseteq K$ also shows that $\Gamma_A$ is pointed. Thus the desired chart is obtained in characteristic zero on the ordinary Nash blowup, and in positive characteristic on its normalization.
\end{enumerate}

This verifies all the required properties of $P$ and the chart indexed by $A$, completing the proof of Theorem~\ref{thm:main}.

Applying Theorem~\ref{thm:main} to $K_1,K_2,K_3$ gives Corollary~\ref{cor:nonnormal-isolated}.

\section{Normal varieties of dimension four with isolated singularities}\label{sec:normal-cycle}

We prove Theorem~\ref{thm:normal-isolated-cycle} by verifying the connecting arrows in~\eqref{eq:normal-chart-paths} and showing that their sources are normal Gorenstein varieties with isolated singularities. The recurring arrows are supplied by Lemma~\ref{lem:recurring-charts}. The sources and intermediate varieties are defined below, in the three characteristic cases.

\subsection{Characteristics other than two and three}\label{subsec:char-not-two-three}

Assume $\operatorname{char}(k)\notin\{2,3\}$.

\paragraph{Source and its properties.}
We define $X_0^{(1)}$ using the lattice
\[
M^{(1)}=\bigl\{(a_1,a_2,a_3,a_4)\in\mathbb Z^4:
 a_1+2a_2+3a_3+5a_4\equiv0\pmod{11}\bigr\}.
\]
Choose the lattice basis
\[
b_1=(11,0,0,0),\quad b_2=(-2,1,0,0),\quad
b_3=(-3,0,1,0),\quad b_4=(-5,0,0,1).
\]
In this basis an exponent $a=(a_1,a_2,a_3,a_4)\in M^{(1)}$ has coordinates
\[
a\longmapsto
\left(\frac{a_1+2a_2+3a_3+5a_4}{11},a_2,a_3,a_4\right)\in\mathbb Z^4.
\]
The cones defining the chart paths in~\eqref{eq:normal-chart-paths} are written in character lattices identified with $\mathbb Z^4$, and their defining matrices list primitive ray generators as columns. A ray is a one-dimensional face of a cone; its primitive generator is the first nonzero lattice point on it.

With respect to the basis $b_1,\ldots,b_4$, the semigroup $M^{(1)}\cap\NN^4$ is identified with $C_0^{(1)}\cap\mathbb Z^4$, where
\begin{equation}\label{eq:source-cone-one}
C_0^{(1)}=\Cone\begin{pmatrix}
1&2&3&5\\
0&11&0&0\\
0&0&11&0\\
0&0&0&11
\end{pmatrix}.
\end{equation}
Set
\[
S_0^{(1)}=C_0^{(1)}\cap\mathbb Z^4,\qquad X_0^{(1)}=\Spec k[S_0^{(1)}].
\]
This is the cyclic quotient $\Aaff_k^4/\mu_{11}$ with weights $(1,2,3,5)$. Here $\mu_r$ denotes the group scheme of $r$-th roots of unity, and the quotient is the spectrum of the invariant ring, including when $\operatorname{char}(k)$ divides $r$.

The source $X_0^{(1)}$ is normal by definition, since $S_0^{(1)}$ is saturated. We verify isolatedness and the Gorenstein property using the dual cone. These lattice arguments are independent of the characteristic and will also apply to $X_0^{(3)}$ in characteristic $3$.

Let $N=\operatorname{Hom}_{\mathbb Z}(\mathbb Z^4,\mathbb Z)$, with the basis dual to the chosen character-lattice basis. The cone $\sigma_0^{(1)}=(C_0^{(1)})^\vee\subset N_{\RR}$ has primitive ray generators
\[
n_1=(11,-2,-3,-5),\quad n_2=(0,1,0,0),\quad
n_3=(0,0,1,0),\quad n_4=(0,0,0,1).
\]
Indeed, if $m_1,\ldots,m_4$ are the displayed ray generators of $C_0^{(1)}$, then $\langle m_j,n_i\rangle=11\delta_{ij}$. The four vectors $n_i$ are primitive and linearly independent, and these pairings identify their cone with $(C_0^{(1)})^\vee$.

By the toric smoothness criterion, $X_0^{(1)}$ is smooth along the orbit corresponding to a face of $\sigma_0^{(1)}$ precisely when that face is regular, meaning that its primitive ray generators extend to a lattice basis of $N$. Since $n_1,\ldots,n_4$ are linearly independent, the facets of $\sigma_0^{(1)}$ are generated by the four triples of these vectors. For a triple of linearly independent lattice vectors, the greatest common divisor (gcd) of the absolute values of the $3\times3$ minors of its column matrix is the index of the lattice they generate in its saturation. The triple therefore extends to a lattice basis exactly when this divisor is one. Here the calculation is
\[
\begin{array}{c|c}
\text{facet generators}&\text{nonzero maximal minors (absolute values)}\\ \hline
n_2,n_3,n_4&1\\
n_1,n_3,n_4&2,11\\
n_1,n_2,n_4&3,11\\
n_1,n_2,n_3&5,11
\end{array}
\]
In each row, the greatest common divisor is one, so all four facets are regular. Every proper face lies in a facet, and every face of a regular cone is regular. Thus $X_0^{(1)}$ is smooth away from its closed torus orbit. On the other hand,
\[
\det\begin{pmatrix}
11&0&0&0\\
-2&1&0&0\\
-3&0&1&0\\
-5&0&0&1
\end{pmatrix}=11.
\]
Hence $\sigma_0^{(1)}$ itself is not regular, so the closed orbit is singular. This orbit is the image of the origin of $\Aaff_k^4$, and it is the unique singular point of $X_0^{(1)}$.

For the Gorenstein property, take the lattice character $m_{\mathrm{can}}=(1,1,1,1)$. It satisfies
\[
\langle m_{\mathrm{can}},n_1\rangle=11-2-3-5=1,
\qquad
\langle m_{\mathrm{can}},n_i\rangle=1\quad(i=2,3,4).
\]
The toric Cartier criterion therefore shows that the canonical divisor of $X_0^{(1)}$ is Cartier. Normal affine toric varieties are Cohen--Macaulay, so $X_0^{(1)}$ is Gorenstein. In the original coordinates, this character corresponds to the invariant monomial $x_1x_2x_3x_4$, whose weight is $1+2+3+5=11$.

\paragraph{Cones and Hilbert bases.}
The path from $X_0^{(1)}$ to $X_{K_1}$ uses the intermediate and recurring cones
\begin{equation}\label{eq:intermediate-cones-one}
\begin{aligned}
C_1^{(1)}&=\Cone\begin{pmatrix}
0&0&1&1&1&2\\
0&1&-2&0&1&-1\\
1&0&1&0&0&0\\
0&0&2&0&2&1
\end{pmatrix},
&\qquad C_2^{(1)}&=\Cone\begin{pmatrix}
-3&-3&0&0&0&1\\
0&1&0&0&1&0\\
4&2&1&2&0&0\\
3&0&0&3&0&0
\end{pmatrix},\\[6pt]
C_3^{(1)}&=\Cone\begin{pmatrix}
0&0&0&1&2\\
0&1&2&0&2\\
1&0&3&0&1\\
0&0&3&0&3
\end{pmatrix},
&\qquad K_1&=\Cone\begin{pmatrix}
0&0&0&1&1&2\\
0&0&1&0&3&3\\
0&1&0&0&-1&-2\\
1&0&0&0&-1&-1
\end{pmatrix}.
\end{aligned}
\end{equation}
The displayed generators of $K_1$ are a reordering of the columns of $B_1$ in~\eqref{eq:recurring-cones}. This is the terminal cone of the first connecting path in~\eqref{eq:normal-chart-paths}.

The cones $C_0^{(1)},C_1^{(1)},C_2^{(1)},C_3^{(1)},K_1$ are pointed, since $x_1+2x_2+2x_3+x_4$ is strictly positive on each of their displayed ray generators.

The source Hilbert basis is
\[
\begin{aligned}
H_0^{(1)}={}&\{(1,b,c,d):b,c,d\in\NN,\ 2b+3c+5d\le11\}\\
&{}\cup\{(2,0,7,0),(2,1,0,4),(2,2,1,3),(2,11,0,0),\\
&\qquad (3,0,1,6),(3,0,11,0),(5,0,0,11)\}.
\end{aligned}
\]
It has $31$ elements, ordered lexicographically as $h_1,\ldots,h_{31}$.

For the cones $C_i^{(a)}$ in these chart paths, put $S_i^{(a)}=C_i^{(a)}\cap\mathbb Z^4$ and let $H_i^{(a)}$ be its Hilbert basis, so $X_i^{(a)}=\Spec k[S_i^{(a)}]$. We use the same lexicographic ordering whenever a source Hilbert basis in this section is indexed as $h_1,\ldots,h_s$, independently of the displayed grouping of its vectors.

The Hilbert bases used in the three connecting paths can be verified by a finite check. For a full-dimensional pointed rational polyhedral cone $C\subset\RR^4$, take its primitive ray generators and the nonzero lattice points in $R[0,1)^4$, where $R$ ranges over matrices whose columns are four linearly independent primitive ray generators of $C$. These vectors generate $C\cap\mathbb Z^4$, by separating the integral and fractional parts of coefficients in the cones $\Cone(R)$, which cover $C$. From this finite set $E$, discard $b$ whenever $b-c\in C$ for some $c\in E\setminus\{b\}$; the remaining vectors form the Hilbert basis. Cone membership is checked by expressing a vector in the real bases given by these matrices and testing whether at least one expression has nonnegative coefficients.

The target column of Table~\ref{tab:chart-generators} lists $H_1^{(1)},H_2^{(1)},H_3^{(1)}$ and the Hilbert basis of $K_1\cap\mathbb Z^4$ in full.

\paragraph{Selected charts and witnesses.}
Table~\ref{tab:ordinary-charts} gives the selected basis $A_i$ and coordinate change $U_i$ for each connecting arrow, together with their determinants and the source Hilbert-basis sizes. For the first arrow, $A_0=[h_3\ h_7\ h_{26}\ h_{29}]$. Every column of $A_i$ belongs to $H_i^{(1)}$, and $\det A_i,\det U_i=\pm1$, so each $A_i$ is a $k$-basis in every characteristic and all coordinate changes are unimodular, that is, invertible over $\mathbb Z$.

\begingroup
\small
\setlength{\tabcolsep}{10pt}
\setlength{\arraycolsep}{3pt}
\renewcommand{\arraystretch}{1.05}
\begin{longtable}{ccc}
\caption{Selected charts in characteristics other than $2$ and $3$.}\label{tab:ordinary-charts}\\
\hline
Arrow and $|H_i^{(1)}|$ & $A_i$ & $U_i$\\\hline
\endfirsthead
\caption[]{Selected charts (continued).}\\
\hline
Arrow and $|H_i^{(1)}|$ & $A_i$ & $U_i$\\\hline
\endhead
\hline\multicolumn{3}{r}{\footnotesize Continued on the next page.}\\
\endfoot
\hline
\endlastfoot
$\begin{gathered}X_0^{(1)}\rightsquigarrow X_1^{(1)}\\ |H_0^{(1)}|=31\end{gathered}$ & $\begin{gathered}\begin{pmatrix}1&1&2&3\\
0&0&1&0\\
0&2&0&1\\
2&1&4&6\end{pmatrix}\\[-2pt]\det A_0=-1\end{gathered}$ & $\begin{gathered}\begin{pmatrix}12&-1&-2&-5\\
-8&0&1&3\\
5&0&-1&-2\\
19&-2&-4&-8\end{pmatrix}\\[-2pt]\det U_0=1\end{gathered}$\\[6pt]
$\begin{gathered}X_1^{(1)}\rightsquigarrow X_2^{(1)}\\ |H_1^{(1)}|=8\end{gathered}$ & $\begin{gathered}\begin{pmatrix}0&1&1&1\\
1&-2&0&1\\
0&1&0&0\\
0&2&1&2\end{pmatrix}\\[-2pt]\det A_1=1\end{gathered}$ & $\begin{gathered}\begin{pmatrix}-2&-1&-3&0\\
0&0&1&0\\
4&2&4&-1\\
3&2&3&-1\end{pmatrix}\\[-2pt]\det U_1=1\end{gathered}$\\[6pt]
$\begin{gathered}X_2^{(1)}\rightsquigarrow X_3^{(1)}\\ |H_2^{(1)}|=12\end{gathered}$ & $\begin{gathered}\begin{pmatrix}0&0&0&1\\
0&0&1&0\\
1&1&0&0\\
0&1&0&0\end{pmatrix}\\[-2pt]\det A_2=-1\end{gathered}$ & $\begin{gathered}\begin{pmatrix}2&3&3&-1\\
2&2&4&-1\\
1&1&3&-1\\
3&3&6&-2\end{pmatrix}\\[-2pt]\det U_2=1\end{gathered}$\\[6pt]
$\begin{gathered}X_3^{(1)}\rightsquigarrow X_{K_1}\\ |H_3^{(1)}|=8\end{gathered}$ & $\begin{gathered}\begin{pmatrix}0&0&1&2\\
0&1&0&2\\
1&1&0&1\\
0&1&0&3\end{pmatrix}\\[-2pt]\det A_3=-1\end{gathered}$ & $\begin{gathered}\begin{pmatrix}0&2&1&-1\\
0&3&0&-1\\
0&-1&0&0\\
1&0&1&-1\end{pmatrix}\\[-2pt]\det U_3=-1\end{gathered}$\\[6pt]
\end{longtable}
\endgroup

In Tables~\ref{tab:ordinary-charts}, \ref{tab:characteristic-three-charts}, and~\ref{tab:characteristic-two-charts}, the selected bases are ordered column matrices, and the coordinate changes act on column vectors. The tabulated vectors and matrices have integer entries; only the determinant tests are performed in $k$.

Table~\ref{tab:chart-generators} gives witnesses for every target Hilbert-basis vector.

\begingroup
\small
\setlength{\tabcolsep}{8pt}
\renewcommand{\arraystretch}{1.1}
\begin{longtable}{cccc}
\caption{Representations of the target Hilbert-basis generators.}\label{tab:chart-generators}\\
\hline
$t$ & $h$ & $q$ & $\delta$\\
\hline
\endfirsthead
\caption[]{Representations of the target Hilbert-basis generators (continued).}\\
\hline
$t$ & $h$ & $q$ & $\delta$\\
\hline
\endhead
\hline
\multicolumn{4}{r}{\footnotesize Continued on the next page.}\\
\endfoot
\hline
\endlastfoot
\multicolumn{4}{c}{$X_0^{(1)}\rightsquigarrow X_1^{(1)}$}\\*
$(0,0,1,0)$ & $(1,3,0,1)$ & $4$ & $2$\\*
$(0,1,0,0)$ & $(1,1,1,1)$ & $4$ & $1$\\
$(1,-2,1,2)$ & $(5,0,0,11)$ & $4$ & $-2$\\
$(1,0,0,0)$ & $(1,1,3,0)$ & $4$ & $1$\\
$(1,1,0,2)$ & $(1,0,1,1)$ & $4$ & $1$\\
$(2,-1,0,1)$ & $(1,1,3,0)$ & $3$ & $-1$\\
$(1,-1,1,1)$ & $(2,2,1,3)$ & $4$ & $1$\\
$(1,0,0,1)$ & $(1,1,1,1)$ & $3$ & $-1$\\
\hline
\multicolumn{4}{c}{$X_1^{(1)}\rightsquigarrow X_2^{(1)}$}\\*
$(-3,0,4,3)$ & $(1,1,0,2)$ & $0$ & ---\\*
$(-3,1,2,0)$ & $(1,-2,1,2)$ & $0$ & ---\\
$(-2,0,3,2)$ & $(1,0,0,1)$ & $0$ & ---\\
$(-2,1,2,1)$ & $(1,-1,1,1)$ & $3$ & $1$\\
$(-1,0,2,1)$ & $(1,0,0,0)$ & $1$ & $1$\\
$(-1,0,2,2)$ & $(0,1,0,0)$ & $0$ & ---\\
$(-1,1,1,0)$ & $(1,-1,1,1)$ & $4$ & $-1$\\
$(0,0,1,0)$ & $(2,-1,0,1)$ & $4$ & $-1$\\
$(0,0,1,1)$ & $(1,0,0,0)$ & $3$ & $2$\\
$(0,0,2,3)$ & $(0,0,1,0)$ & $2$ & $1$\\
$(0,1,0,0)$ & $(0,0,1,0)$ & $4$ & $-1$\\
$(1,0,0,0)$ & $(1,0,0,0)$ & $4$ & $-1$\\
\hline
\multicolumn{4}{c}{$X_2^{(1)}\rightsquigarrow X_3^{(1)}$}\\*
$(0,0,1,0)$ & $(-3,1,2,0)$ & $1$ & $-2$\\*
$(0,1,0,0)$ & $(0,0,2,3)$ & $1$ & $1$\\
$(0,1,1,1)$ & $(-1,0,2,1)$ & $1$ & $-1$\\
$(0,2,3,3)$ & $(-3,1,2,0)$ & $3$ & $-1$\\
$(1,0,0,0)$ & $(-1,1,1,0)$ & $1$ & $-1$\\
$(1,2,1,2)$ & $(0,0,2,3)$ & $2$ & $-3$\\
$(1,2,2,3)$ & $(-1,0,2,1)$ & $2$ & $-1$\\
$(2,2,1,3)$ & $(1,0,0,0)$ & $0$ & ---\\
\hline
\multicolumn{4}{c}{$X_3^{(1)}\rightsquigarrow X_{K_1}$}\\*
$(0,0,0,1)$ & $(1,0,0,0)$ & $0$ & ---\\*
$(0,0,1,0)$ & $(0,1,0,0)$ & $4$ & $1$\\
$(0,1,0,0)$ & $(0,1,0,0)$ & $2$ & $-3$\\
$(1,0,0,0)$ & $(1,2,2,3)$ & $4$ & $-1$\\
$(1,2,-1,0)$ & $(1,2,1,2)$ & $2$ & $-2$\\
$(1,3,-1,-1)$ & $(0,1,0,0)$ & $1$ & $2$\\
$(2,3,-2,-1)$ & $(1,2,2,3)$ & $1$ & $-1$\\
\end{longtable}
\endgroup

For each block of Table~\ref{tab:chart-generators}, let $C,C'$ be the source and target cones, and let $H,H'$ be the Hilbert bases of $C\cap\mathbb Z^4$ and $C'\cap\mathbb Z^4$, respectively. Put $S'=C'\cap\mathbb Z^4$ and take $A,U$ from the corresponding row of Table~\ref{tab:ordinary-charts}. The vectors $t\in H'$ and $h\in H$ satisfy
\[
t=\begin{cases}
Uh,&q=0,\\
U(h-a_q),&q>0,
\end{cases}
\]
where $a_q$ is column $q$ of $A$. Thus $q=0$ denotes an original generator, while $q\in\{1,2,3,4\}$ specifies the selected column used in the exchange. In the latter case, $h\in H\setminus A$ and $\delta$ is the determinant obtained by replacing column $q$ of $A$ by $h$; the exchange is valid when $\delta\ne0$ in $k$. A dash in the $\delta$ column marks an original generator. We retain the meanings of $q$ and $\delta$ in Tables~\ref{tab:characteristic-three-witnesses} and~\ref{tab:characteristic-two-witnesses}.

Set $p=\operatorname{char}(k)$ and let $G_{A,p}$ consist of $H$ and all valid exchange differences $h-a$ with $h\in H\setminus A$ and $a\in A$. Write $\Gamma_{A,p}=\langle G_{A,p}\rangle_{\NN}$; the subscript $0$ denotes the characteristic-zero determinant test. These are the generators and semigroup of the selected chart in Lemma~\ref{lem:chart}.

Every exchange in Table~\ref{tab:chart-generators} has determinant $\pm1$, $\pm2$, or $\pm3$, and therefore survives when $p\notin\{2,3\}$; original generators always survive. Thus the witnesses give $H'\subseteq UG_{A,p}$ for each connecting arrow.

The inclusion $UG_{A,p}\subseteq S'$ is verified by enumerating $G_{A,p}$, applying $U$, and checking membership in $C'$ using the cone-membership test described after $H_0^{(1)}$; the transformed vectors are integral. For the four arrows in Table~\ref{tab:ordinary-charts}, it suffices to check $G_{A,0}$, since $G_{A,p}\subseteq G_{A,0}$.

\paragraph{Conclusion.}
The witnesses and this containment check give
\[
H'\subseteq UG_{A,p}\subseteq S',
\qquad U\Gamma_{A,p}=S'.
\]
Since $C'$ is pointed, so is $\Gamma_{A,p}$. Lemma~\ref{lem:chart} identifies the selected ordinary chart with $\Spec k[S']$, proving all four connecting arrows. These charts are normal, so the same arrows also hold for normalized Nash blowups. The recurrence of $X_{K_1}$ follows from Lemma~\ref{lem:recurring-charts}.

\subsection{Characteristic three}\label{subsec:characteristic-three}

Assume $\operatorname{char}(k)=3$.

\paragraph{Source and its properties.}
Since $C_0^{(3)}=C_0^{(1)}$, the characteristic-independent argument in Section~\ref{subsec:char-not-two-three} proves that $X_0^{(3)}$ is normal and Gorenstein with an isolated singularity.

\paragraph{Cones and Hilbert bases.}
The connecting path uses the cones
\begin{equation}\label{eq:connecting-cones-three}
C_2^{(3)}=\Cone\begin{pmatrix}-3&0&0&0&1&1\\3&0&0&1&-1&0\\1&0&1&0&0&0\\1&1&0&0&1&0\end{pmatrix},\qquad
K'_3=\Cone\begin{pmatrix}0&0&0&1&1\\0&1&2&0&1\\0&0&3&0&3\\1&0&0&-1&0\end{pmatrix}.
\end{equation}
The cone $K'_3$ is denoted by $\omega_{3,1}$ in~\cite[Section~2]{counterexamples}; thus $X_{K'_3}$ is the intermediate chart in Lemma~\ref{lem:recurring-charts}.

The cones $C_1^{(3)},C_2^{(3)},K'_3$ are pointed, since $2x_1+2x_2+x_3+x_4$ is positive on each of their displayed ray generators. Since $C_i^{(3)}=C_i^{(1)}$ for $i=0,1$, we have $H_i^{(3)}=H_i^{(1)}$ for these indices. The finite enumeration described in Section~\ref{subsec:char-not-two-three} gives the remaining source Hilbert basis
\begingroup\small
\[
H_2^{(3)}=\left\{\begin{gathered}
(-3,3,1,1),\ (-2,2,1,1),\ (-1,1,1,1),\ (0,0,0,1),\\
(0,0,1,0),\ (0,1,0,0),\ (1,-1,0,1),\ (1,0,0,0)
\end{gathered}\right\}.
\]
\endgroup

\paragraph{Selected charts and witnesses.}
For the first arrow $X_0^{(3)}\overset{n}{\rightsquigarrow}X_1^{(3)}$, use $A_0,U_0$ from Table~\ref{tab:ordinary-charts}. Its witnesses in the first block of Table~\ref{tab:chart-generators} have replacement determinants $\pm1$ or $\pm2$, so they survive in characteristic $3$ and give $H_1^{(3)}\subseteq U_0G_{A_0,3}$.

Table~\ref{tab:characteristic-three-charts} records the two remaining chart choices in the format of Table~\ref{tab:ordinary-charts}. Again, the columns of $A$ lie in the source Hilbert basis and $\det A,\det U=\pm1$.

\begingroup
\small
\setlength{\tabcolsep}{10pt}
\setlength{\arraycolsep}{3pt}
\renewcommand{\arraystretch}{1.05}
\begin{longtable}{ccc}
\caption{Selected charts in characteristic $3$.}\label{tab:characteristic-three-charts}\\
\hline
Arrow and $|H|$ & $A$ & $U$\\\hline
\endfirsthead
\caption[]{Selected charts in characteristic $3$ (continued).}\\
\hline
Arrow and $|H|$ & $A$ & $U$\\\hline
\endhead
\hline\multicolumn{3}{r}{\footnotesize Continued on the next page.}\\
\endfoot
\hline
\endlastfoot
$\begin{gathered}X_1^{(3)}\overset{n}{\rightsquigarrow}X_2^{(3)}\\ |H|=8\end{gathered}$ & $\begin{gathered}\begin{pmatrix}0&1&1&1\\1&-2&0&0\\0&1&0&0\\0&2&0&1\end{pmatrix}\\[-2pt]\det A=1\end{gathered}$ & $\begin{gathered}\begin{pmatrix}0&-1&0&-1\\0&1&1&1\\1&1&1&0\\2&1&2&0\end{pmatrix}\\[-2pt]\det U=-1\end{gathered}$\\[6pt]
$\begin{gathered}X_2^{(3)}\overset{n}{\rightsquigarrow}X_{K'_3}\\ |H|=8\end{gathered}$ & $\begin{gathered}\begin{pmatrix}-3&0&0&1\\3&0&1&0\\1&1&0&0\\1&0&0&0\end{pmatrix}\\[-2pt]\det A=1\end{gathered}$ & $\begin{gathered}\begin{pmatrix}1&1&0&1\\1&0&1&2\\1&0&0&3\\-1&-1&0&0\end{pmatrix}\\[-2pt]\det U=1\end{gathered}$\\[6pt]
\end{longtable}
\endgroup

Table~\ref{tab:characteristic-three-witnesses} gives witnesses for every primitive ray generator of the target cones.

\begingroup
\small
\setlength{\tabcolsep}{8pt}
\setlength{\arraycolsep}{3pt}
\renewcommand{\arraystretch}{1.1}
\begin{longtable}{cccc}
\caption{Representations of the primitive target ray generators in characteristic $3$.}\label{tab:characteristic-three-witnesses}\\
\hline
$t$ & $h$ & $q$ & $\delta$\\\hline
\endfirsthead
\caption[]{Representations of the primitive target ray generators in characteristic $3$ (continued).}\\
\hline
$t$ & $h$ & $q$ & $\delta$\\\hline
\endhead
\hline\multicolumn{4}{r}{\footnotesize Continued on the next page.}\\
\endfoot
\hline
\endlastfoot
\multicolumn{4}{c}{$X_1^{(3)}\overset{n}{\rightsquigarrow}X_2^{(3)}$}\\*
$(-3,3,1,1)$ & $(1,1,0,2)$ & $3$ & $-1$\\*
$(0,0,0,1)$ & $(2,-1,0,1)$ & $3$ & $1$\\
$(0,0,1,0)$ & $(0,0,1,0)$ & $2$ & $1$\\
$(0,1,0,0)$ & $(0,0,1,0)$ & $3$ & $1$\\
$(1,-1,0,1)$ & $(2,-1,0,1)$ & $4$ & $1$\\
$(1,0,0,0)$ & $(0,0,1,0)$ & $4$ & $-2$\\
\hline
\multicolumn{4}{c}{$X_2^{(3)}\overset{n}{\rightsquigarrow}X_{K'_3}$}\\*
$(0,0,0,1)$ & $(-2,2,1,1)$ & $4$ & $1$\\*
$(0,1,0,0)$ & $(0,0,1,0)$ & $0$ & ---\\
$(0,2,3,0)$ & $(0,0,0,1)$ & $1$ & $1$\\
$(1,0,0,-1)$ & $(0,1,0,0)$ & $0$ & ---\\
$(1,1,3,0)$ & $(0,0,0,1)$ & $2$ & $-1$\\
\end{longtable}
\endgroup

Use the column conventions after Table~\ref{tab:chart-generators}, with $A,U$ from Table~\ref{tab:characteristic-three-charts}; here $t$ ranges over the primitive ray generators $\rho_1,\ldots,\rho_m$ of the target cone $C'$. The exchange determinants are nonzero modulo $3$, so the witnesses give $\{\rho_1,\ldots,\rho_m\}\subseteq UG_{A,3}$ for both arrows.

\paragraph{Conclusion.}
The finite check following Table~\ref{tab:chart-generators} gives $UG_{A,3}\subseteq C'$ for each connecting arrow. Together with the witnesses, this identifies the cone of the transformed chart semigroup with $C'$. The semigroup $\Gamma_{A,3}$ is pointed because $C'$ is pointed, and it generates $\mathbb Z^4$ because it contains the unimodular basis $A$. Consequently,
\[
U\Cone(\Gamma_{A,3})=C',
\qquad U\overline{\Gamma_{A,3}}=C'\cap\mathbb Z^4.
\]
Lemma~\ref{lem:chart} identifies the normalization of the selected chart with $\Spec k[C'\cap\mathbb Z^4]$. This proves all three connecting arrows, reaching $X_{K'_3}$ from $X_0^{(3)}$. The last two arrows in the characteristic-$3$ path of~\eqref{eq:normal-chart-paths} follow from Lemma~\ref{lem:recurring-charts}.

\subsection{Characteristic two}\label{subsec:characteristic-two}

Assume $\operatorname{char}(k)=2$.

\paragraph{Source and its properties.}
Let
\[
M^{(2)}=\{a\in\mathbb Z^4:a_1+2a_2+7a_3+8a_4\equiv0\pmod9\}.
\]
In the basis $(9,0,0,0),(-2,1,0,0),(-7,0,1,0),(-8,0,0,1)$, the semigroup $M^{(2)}\cap\NN^4$ has cone
\begin{equation}\label{eq:source-cone-two}
C_0^{(2)}=\Cone\begin{pmatrix}1&2&7&8\\0&9&0&0\\0&0&9&0\\0&0&0&9\end{pmatrix}.
\end{equation}
Put $X_0^{(2)}=\Spec k[S_0^{(2)}]$. This is the cyclic quotient $\Aaff_k^4/\mu_9$ with weights $(1,2,7,8)$.

The source $X_0^{(2)}$ is normal by definition, since $S_0^{(2)}$ is saturated. The primitive ray generators of the dual cone of $C_0^{(2)}$ are $(9,-2,-7,-8),e_2,e_3,e_4$, where $e_i$ denotes the $i$-th standard basis vector of the dual lattice. The nonzero facet minors have absolute values $(9,2)$, $(9,7)$, $(9,8)$, and $1$, so every proper face is regular; the full determinant has absolute value $9$. Thus $X_0^{(2)}$ has an isolated singularity. The character $(2,1,1,1)$ pairs to $1$ with these four generators, so the Cartier and Cohen--Macaulay arguments of Section~\ref{subsec:char-not-two-three} prove that $X_0^{(2)}$ is Gorenstein.

\paragraph{Cones and Hilbert bases.}
The connecting path uses
\begin{equation}\label{eq:intermediate-cone-two}
C_1^{(2)}=\Cone\begin{pmatrix}1&0&0&0&5&-1&1\\0&1&0&0&-2&1&-2\\0&0&1&0&5&1&5\\0&0&0&1&0&0&4\end{pmatrix},
\qquad X_1^{(2)}=\Spec k[S_1^{(2)}].
\end{equation}
The terminal cone $K_2$ is given in~\eqref{eq:recurring-cones}. The cones $C_0^{(2)},C_1^{(2)},K_2$ are pointed, since $2x_1+2x_2+x_3+x_4$ is positive on each of their displayed ray generators. The finite enumeration described in Section~\ref{subsec:char-not-two-three} gives the source Hilbert bases
\begingroup\small
\[
H_0^{(2)}=\left\{\begin{gathered}
(1,0,0,0),\ (1,0,0,1),\ (1,0,1,0),\ (1,1,0,0),\\
(1,1,1,0),\ (1,2,0,0),\ (1,3,0,0),\ (1,4,0,0),\\
(2,1,0,2),\ (2,5,0,1),\ (2,9,0,0),\ (4,0,4,1),\\
(4,0,5,0),\ (5,0,3,3),\ (6,0,2,5),\ (7,0,1,7),\\
(7,0,9,0),\ (8,0,0,9)
\end{gathered}\right\}.
\]
\[
H_1^{(2)}=\left\{\begin{gathered}
(-1,1,1,0),\ (0,0,0,1),\ (0,0,1,0),\ (0,1,0,0),\\
(1,-2,5,4),\ (1,-1,3,2),\ (1,0,0,0),\ (2,-2,5,3),\\
(2,-1,3,1),\ (3,-2,5,2),\ (3,-1,3,0),\ (4,-2,5,1),\\
(5,-2,5,0)
\end{gathered}\right\}.
\]
\endgroup

\paragraph{Selected charts and witnesses.}
Table~\ref{tab:characteristic-two-charts} records the two connecting chart choices in the format of Table~\ref{tab:ordinary-charts}. Again, the columns of $A$ lie in the source Hilbert basis and $\det A,\det U=\pm1$.

\begingroup
\small
\setlength{\tabcolsep}{10pt}
\setlength{\arraycolsep}{3pt}
\renewcommand{\arraystretch}{1.05}
\begin{longtable}{ccc}
\caption{Selected charts in characteristic $2$.}\label{tab:characteristic-two-charts}\\
\hline
Arrow and $|H|$ & $A$ & $U$\\\hline
\endfirsthead
\caption[]{Selected charts in characteristic $2$ (continued).}\\
\hline
Arrow and $|H|$ & $A$ & $U$\\\hline
\endhead
\hline\multicolumn{3}{r}{\footnotesize Continued on the next page.}\\
\endfoot
\hline
\endlastfoot
$\begin{gathered}X_0^{(2)}\overset{n}{\rightsquigarrow}X_1^{(2)}\\ |H|=18\end{gathered}$ & $\begin{gathered}\begin{pmatrix}1&1&2&1\\1&0&1&4\\1&0&0&0\\0&1&2&0\end{pmatrix}\\[-2pt]\det A=-1\end{gathered}$ & $\begin{gathered}\begin{pmatrix}4&-1&-3&-3\\1&0&-1&-1\\2&0&0&-1\\0&0&1&0\end{pmatrix}\\[-2pt]\det U=-1\end{gathered}$\\[6pt]
$\begin{gathered}X_1^{(2)}\overset{n}{\rightsquigarrow}X_{K_2}\\ |H|=13\end{gathered}$ & $\begin{gathered}\begin{pmatrix}1&0&2&1\\-2&0&-2&-1\\5&1&5&3\\4&0&3&2\end{pmatrix}\\[-2pt]\det A=-1\end{gathered}$ & $\begin{gathered}\begin{pmatrix}2&3&0&1\\3&10&2&2\\4&11&2&2\\1&5&1&1\end{pmatrix}\\[-2pt]\det U=-1\end{gathered}$\\[6pt]
\end{longtable}
\endgroup

Table~\ref{tab:characteristic-two-witnesses} gives witnesses for every primitive ray generator of the target cones.

\begingroup
\small
\setlength{\tabcolsep}{8pt}
\setlength{\arraycolsep}{3pt}
\renewcommand{\arraystretch}{1.1}
\begin{longtable}{cccc}
\caption{Representations of the primitive target ray generators in characteristic $2$.}\label{tab:characteristic-two-witnesses}\\
\hline
$t$ & $h$ & $q$ & $\delta$\\\hline
\endfirsthead
\caption[]{Representations of the primitive target ray generators in characteristic $2$ (continued).}\\
\hline
$t$ & $h$ & $q$ & $\delta$\\\hline
\endhead
\hline\multicolumn{4}{r}{\footnotesize Continued on the next page.}\\
\endfoot
\hline
\endlastfoot
\multicolumn{4}{c}{$X_0^{(2)}\overset{n}{\rightsquigarrow}X_1^{(2)}$}\\*
$(1,0,0,0)$ & $(1,0,1,0)$ & $1$ & $-1$\\*
$(0,1,0,0)$ & $(1,3,0,0)$ & $3$ & $1$\\
$(0,0,1,0)$ & $(2,5,0,1)$ & $4$ & $-1$\\
$(0,0,0,1)$ & $(1,0,1,0)$ & $3$ & $1$\\
$(5,-2,5,0)$ & $(8,0,0,9)$ & $4$ & $1$\\
$(-1,1,1,0)$ & $(2,5,0,1)$ & $3$ & $-1$\\
$(1,-2,5,4)$ & $(4,0,4,1)$ & $4$ & $1$\\
\hline
\multicolumn{4}{c}{$X_1^{(2)}\overset{n}{\rightsquigarrow}X_{K_2}$}\\*
$(1,0,0,0)$ & $(0,0,0,1)$ & $2$ & $1$\\*
$(0,1,0,0)$ & $(1,-2,5,4)$ & $0$ & ---\\
$(1,1,2,0)$ & $(2,-1,3,1)$ & $4$ & $-1$\\
$(0,0,0,1)$ & $(0,0,0,1)$ & $3$ & $1$\\
$(0,2,2,1)$ & $(0,0,1,0)$ & $0$ & ---\\
\end{longtable}
\endgroup

Use the target-cone notation and column conventions after Table~\ref{tab:characteristic-three-witnesses}, with $A,U$ from Table~\ref{tab:characteristic-two-charts}. The exchange determinants are $\pm1$, so the witnesses give $\{\rho_1,\ldots,\rho_m\}\subseteq UG_{A,2}$ for both arrows.

\paragraph{Conclusion.}
The containment check after Table~\ref{tab:chart-generators} gives $UG_{A,2}\subseteq C'$ for each arrow. Since the target cones are pointed and each $\Gamma_{A,2}$ contains the unimodular basis $A$, the normalized-chart argument of Section~\ref{subsec:characteristic-three} applies, proving both connecting arrows to $X_{K_2}$. Its recurrence follows from Lemma~\ref{lem:recurring-charts}.

Normalization is essential for the selected arrow $X_1^{(2)}\overset{n}{\rightsquigarrow}X_{K_2}$: the set $UG_{A,2}$ misses the Hilbert-basis elements $(0,1,1,1)$ and $(1,1,1,0)$. Their indecomposability, recorded in the Hilbert basis of $K_2\cap\mathbb Z^4$ in Appendix~\ref{app:ordinary-refinement}, prevents them from being sums of transformed generators. Thus the corresponding ordinary chart is not normal.

\subsection{Conclusion and computational remarks}

The connecting paths, together with the recurrences in Lemma~\ref{lem:recurring-charts}, prove Theorem~\ref{thm:normal-isolated-cycle}. Take $X_0^{(2)}$ in characteristic $2$, $X_0^{(3)}$ in characteristic $3$, and $X_0^{(1)}$ otherwise. By Sections~\ref{subsec:char-not-two-three} and~\ref{subsec:characteristic-two}, these sources are normal Gorenstein affine toric varieties of dimension $4$ with an isolated singularity. Since Nash blowup and normalization commute with restriction to open subsets, the paths and recurrences produce singular charts at arbitrarily late stages of normalized Nash iteration. No earlier iterate can be smooth: both operations are isomorphisms on smooth varieties, so a smooth iterate would force every subsequent iterate to be smooth. If $\operatorname{char}(k)\notin\{2,3\}$, the selected ordinary charts are also normal, and the same argument applies to ordinary Nash iteration.

\begin{remark}
The normal sources were found among cyclic quotients $\Aaff_k^4/\mu_r$ with weights $(1,a,b,c)$ satisfying
\[
r\ge2,\qquad 1\le a\le b\le c<r,\qquad \gcd(r,abc)=1.
\]
The two sources used here have weights $(1,2,3,5)$ and $(1,2,7,8)$, with orders $11$ and $9$, respectively.

This family offers several advantages for the search:
\begin{enumerate}
\item The quotients are normal, and the coprimality condition ensures an isolated singularity.
\item Each defining cone has four primitive ray generators and is determined by the order and weights.
\item For each $r$ there are finitely many weight tuples. Coordinate permutations and simultaneous multiplication of the weights modulo $r$ by an integer coprime to $r$ reduce repeated candidates.
\item The Gorenstein condition is the congruence $1+a+b+c\equiv0\pmod r$.
\end{enumerate}

The SageMath programs by Castillo et al., described in~\cite[Section~3]{computational}, provided valuable computational support for the search.
\end{remark}

\raggedbottom
\section{Conclusions and remaining questions}

The two results are complementary. Theorem~\ref{thm:main} gives a uniform construction from prescribed cones in arbitrary dimension. Its sources are nonnormal, with isolated singularities and smooth normalization. Theorem~\ref{thm:normal-isolated-cycle} gives normal Gorenstein affine toric examples in dimension $4$ through specific chart paths. Both types of isolated counterexamples remain singular under iterated normalized Nash blowups in every characteristic. Their ordinary iterates are also singular in characteristic zero for the nonnormal examples and in characteristics other than $2$ and $3$ for the normal examples. Appendix~\ref{app:ordinary-refinement} gives an ordinary-chart refinement in selected positive characteristics.

These results answer the open problem of Castillo, Duarte, Leyton-\'Alvarez, and Liendo for isolated singularities in dimension $4$ and suggest further questions about the mechanism behind the examples. The behavior of ordinary Nash iteration for the displayed normal sources in characteristics $2$ and $3$ remains to be determined.

In characteristic zero, which geometric or combinatorial properties of a cone and its lattice-point semigroup force iterated normalized Nash blowups to reach a singular recurring chart? Can such properties be recognized before computing the full chart path and used to construct families of normal affine toric varieties of dimension $4$ with isolated singularities for which neither ordinary nor normalized Nash blowup iteration gives a resolution?

In lower dimensions, the questions of whether iterated ordinary Nash blowups resolve all surfaces in characteristic zero and whether iterated normalized Nash blowups resolve all normal varieties of dimension $3$ remain open; see~\cite[Sections~4.1 and~4.2]{computational}.

\appendix
\section{Hilbert bases and ordinary-chart refinement}\label{app:ordinary-refinement}

Consider the cones $K_1,K_2,K_3$ defined in~\eqref{eq:recurring-cones}. Set $M=\mathbb Z^4$. The ordered Hilbert bases $\mathcal H_i=[h_{i,1},\ldots,h_{i,s_i}]$ of $K_i\cap M$, for $i=1,2,3$, are
\begin{small}
\begin{align*}
\mathcal H_1=[{}&(1,0,0,0),(0,1,0,0),(0,0,1,0),(0,0,0,1),\\
&(2,3,-2,-1),(1,3,-1,-1),(1,2,-1,0)],\\
\mathcal H_2=[{}&(1,0,0,0),(0,1,0,0),(1,1,2,0),(0,0,0,1),\\
&(0,2,2,1),(1,1,1,0),(0,1,1,1)],\\
\mathcal H_3=[{}&(1,0,0,0),(0,1,0,0),(1,2,3,0),(0,0,0,1),\\
&(0,3,3,1),(0,2,2,1),(1,1,1,0),(1,2,2,0),(0,1,1,1)].
\end{align*}
\end{small}

We now describe a refinement of Theorem~\ref{thm:main} that realizes $X_{K_1}$ in characteristic $p\ge5$ and $X_{K_3}$ in characteristic $3$ as ordinary Nash charts.

In the construction of Section~\ref{sec:construction}, write each Hilbert-basis element as $h_j=\sum_i h_j^i v_i$ in the chosen interior lattice basis. Instead of the shift and congruence condition in \eqref{eq:shift} and \eqref{eq:congruence}, use
\[
u_j=u_0+h_j,\qquad t=u_0^d\bmod p.
\]
In the following residue conditions, we regard the coordinates $h_j^d$ as elements of $\mathbb F_p$. The selected set $A$ is a $k$-basis if $t\ne0$, and replacing $u_0$ by $u_j$ is valid if $t+h_j^d\ne0$. Thus the construction gives the ordinary chart $\Spec k[S_K]$ whenever
\[
\{0,h_1^d,\ldots,h_s^d\}\ne\mathbb F_p.
\]
Choose $t$ outside the negatives of these residues and take $u_0$ sufficiently large in the prescribed residue classes. The exchanges replacing $u_0$ by $u_j$ give $u_j-u_0=h_j$, which generate $S_K$.

For $K_1$ in characteristic $p\ge5$, the primitive functional $\lambda(x)=x_1-x_2$ gives the last coordinates of its Hilbert basis modulo $p$ as
\[
(1,-1,0,0,-1,-2,-1).
\]
The forbidden residues for $t$ are therefore $\{0,-1,1,2\}$ modulo $p$, so $t=3$ is allowed for every $p\ge5$.

For $K_3$ in characteristic $3$, use $\lambda(x)=x_1-3x_2$, whose Hilbert-basis values modulo $3$ are
\[
(1,0,1,0,0,0,1,1,0).
\]
The forbidden residues for $t$ are $\{0,2\}$ modulo $3$, so we may take $t=1$.

Both functionals are primitive and take both signs on their cones. For either $\lambda$, its kernel therefore meets the interior of the cone; choose a primitive lattice vector there and extend it to a basis of the kernel lattice, moving the other basis vectors into the interior by adding large multiples of the first. A lattice vector $w$ with $\lambda(w)=1$, adjusted by a large multiple of the same interior vector, completes this to an interior basis of $M$ realizing $\lambda$ as its last-coordinate functional. Thus the residue choices are compatible with the construction, and the prescribed charts occur on the ordinary Nash blowup.

For $K_2$ in characteristic $2$, this residue criterion cannot work for any last-coordinate functional: the reductions of $h_{2,1},h_{2,2},h_{2,4},h_{2,6}$ span $M/2M$. Hence every nonzero functional takes value $1$ on some Hilbert-basis element, and the residues together with zero exhaust $\mathbb F_2$. This obstruction applies only to the refinement described here; other ordinary-chart constructions remain possible.

\end{document}